\documentclass[runningheads]{llncs}
\usepackage{subfig}
\usepackage{graphicx}
\usepackage{amsmath,amssymb, stmaryrd}
\usepackage{tikz}
\usepackage{bm}
\usepackage{hyperref}

\begin{document}
\title{Correlated Boolean Operators for Uncertainty Logic}
%
%
\author{Enrique Miralles-Dolz\inst{1,2}\and 
Ander Gray\inst{1,2} \and
Edoardo Patelli\inst{3} \and
Scott Ferson\inst{1}}
\authorrunning{Miralles-Dolz, Gray et al.}
%
\institute{Institute for Risk and Uncertainty, University of Liverpool
\email{\{enmidol,akgray,ferson\}@liverpool.ac.uk}\and
United Kingdom Atomic Energy Authority
\and
Centre for Intelligent Infrastructure, University of Strathclyde\\
\email{edoardo.patelli@strath.ac.uk}}
\maketitle              
\begin{abstract}

We present a correlated \textit{and} gate which may be used to propagate uncertainty and dependence through Boolean functions, since any Boolean function may be expressed as a combination of \textit{and} and \textit{not} operations. 
We argue that the \textit{and} gate is a bivariate copula family, which has the interpretation of constructing bivariate Bernoulli random variables following a given Pearson correlation coefficient and marginal probabilities. We show how this copula family may be used to propagate uncertainty in the form of probabilities of events, probability intervals, and probability boxes, with only partial or no knowledge of the dependency between events, expressed as an interval for the correlation coefficient. These results generalise previous results by Fréchet on the conjunction of two events with unknown dependencies. We show an application propagating uncertainty through a fault tree for a pressure tank. This paper comes with an open-source Julia library for performing uncertainty logic.

\end{abstract}
\section{Introduction}
\footnotetext[0]{This version of the article has been accepted for publication, after peer review and is subject to Springer Nature’s \href{https://www.springernature.com/gp/open-research/policies/accepted-manuscript-terms}{AM terms of use}, but is not the Version of Record and does not reflect post-acceptance improvements, or any corrections. The Version of Record is available online at: \href{https://doi.org/10.1007/978-3-031-08971-8_64}{https://doi.org/10.1007/978-3-031-08971-8\_64}}
The logical conjunction ($\land$) is a function $f:\{0,1\}^2 \rightarrow{\{0,1\}}$ that returns a value of 1 if and only if
both of the inputs are 1.
%
\noindent The logical values assigned to the inputs represent the truth value of certain propositions or events, and in Boolean algebra it
is required absolute certainty about these truth values (i.e. they are either \textit{true} (1) or \textit{false} (0)). That is, variables take the form $p \in \{0, 1\}$. However, this requirement is often too restrictive to be used in practical applications where truth values have some degree of uncertainty, and therefore an extension of classical Booleans to uncertain Booleans would be desirable. We define uncertain Booleans as probabilities $p \in [0,1]$ defining a precise Bernoulli distribution (an event), interval probabilities $p \subseteq [0,1]$ defining a set of Bernoulli distributions (a credal set), and probability boxes with $\text{range}(p) \subseteq [0,1]$. This paper presents a method to perform uncertainty logic with these structures.

Note that it is possible to build any binary Boolean operation from two primitive operations, for example from \textit{and} and \textit{not}. Table \ref{Tabel1} shows the 16 ($2^4$) possible binary and unary Boolean operations, written in terms of two primitives operations $\land$ and $\neg$. Note that some of these operations (e.g. Identity and Zero) are trivial, and with the other operations being written in terms of previously derived operations for brevity. Since all of these operations may be written in terms of an $\land$ and a $\neg$ operation, and since the \textit{not} operation is $not(A) = \neg A = 1 - A$, it follows that we therefore only need to describe a correlated \textit{and} operator, and all other binary operators, and further more complicated Boolean functions, may be calculated in terms of these two operators.
\begin{table}[!h]
\centering
\caption{Summary of the 16 possible binary and unary Boolean operations written in terms of two primitive operations $\land$ and $\neg$.}
{
\vspace{2em}
\begin{tabular}{ |c|c|c|c|c| }
\hline
    $A$ & $1 1 0 0$  & Adopted Name & Adopted Symbol & Expansion\\ 
    $B$ & $1 0 1 0$  &  &  & \\ 
 \hline
 $w_{1}$ & $1 1 0 0$ & Identity & $A$ & trivial \\ 
 $w_{2}$ & $1 0 1 0$ & Identity & $B$ & trivial \\ 
 $w_{3}$ & $0 0 0 0$ & Zero & 0 & trivial \\ 
 $w_{4}$ & $1 1 1 1$ & One & 1 & trivial \\ 
 $w_{5} $ & $1 0 0 0$ & And & $A\land B$ & primitive \\ 
 $w_{6} $ & $0 0 1 1$ & Not & $\neg A$ & $1 - A$ \\ 
 $w_{7} $ & $0 1 0 1$ & Not & $\neg B$ & $1 - B$ \\ 
 $w_{8} $ & $1 1 1 0$ & Or & $A$ $\lor$ $B$ & $\neg((\neg A) \land (\neg B))$ \\ 
 $w_{9} $ & $0 1 1 1$ & Nand & $A$ nand $B$ & $\neg(A \land B)$ \\ 
 $w_{10}$ & $0 0 0 1$ & Nor & $A$ nor $B$ & $\neg(A$ $\lor$ $B)$ \\ 
 $w_{11}$& $0 1 1 0$ & Exclusive Or  & $A$ xor $B$ & $(A$ $\lor$ $B)\land (A$ nand $B)$ \\ 
 $w_{12}$ & $1 0 0 1$ & Equivalence  & $A\equiv B$ &  $\neg (A$ xor $B)$ \\ 
 $w_{13}$& $1 0 1 1$ & Implication & $A\implies B$ &  $\neg A$ $\lor$ $B$ \\ 
 $w_{14}$ & $1 1 0 1$ & Implication &  $B \implies A$ &  $A$ $\lor$ $\neg B$ \\ 
 $w_{15}$&  $0 1 0 0$ &  Inhibition & $A \Longmapsto B$ &  $\neg (A \implies B)$ \\ 
 $w_{16}$ & $0 0 1 0$ &  Inhibition & $B \Longmapsto A$ & $\neg (B \implies A)$ \\ 
 \hline
\end{tabular}}
\label{Tabel1}
\end{table}

When events $A$ and $B$ are independent, their logical conjunction is calculated as $\mathbb{P}(A \land B) = \mathbb{P}(A)\mathbb{P}(B)$.
However, the assumption of independence has a significant consequence quantitatively, as shown in \cite{ferson2015dependence}. Therefore, after extending the mathematical structures to describe events $A$ and $B$ from classical Booleans to uncertain Booleans, the next desirable extension would be on this assumption of independence between their probabilities. For example, consider the two following random bit-vectors each with the same marginal probabilities $P(A) = P(B) = 0.5$ (a sequence of two fair coin tosses)
\begin{align}
\label{eq:rand_vec}
 A & = \{ 0,  1,  1,  1,  1,  0,  0,  1,  1,  1,  0,  0,  1\}, \notag \\
 B & = \{0,  0,  0,  1,  1, 1, 1,  1,  1,  1,  0,  0,  1\} .
\end{align}
\noindent Although the individual coin tosses is fair (the sample mean for these 12 tosses is $\sim 0.6$ for each vector), the vectors are corrected, that is, the outcome of one throw can influence the other, with $\rho_{AB} = 0.3$ (sample correlation is $\sim 0.35$).

In \cite{lucas1995default}, it is derived a model for the conjunction employing the Pearson correlation coefficient to capture dependence (referred to as the Lucas model). Unlike for continuous distributions, two marginals and a correlation coefficient is sufficient to completely define a bivariate Bernoulli random variable \cite{joe1997multivariate}.
The Lucas model is defined as
\begin{equation}
    \mathbb{P}(A \land B) = \mathbb{P}(A)\mathbb{P}(B) + \rho_{AB}\sqrt{\mathbb{P}(A)\mathbb{P}(\neg A)\mathbb{P}(B)\mathbb{P}(\neg B)} \;,
\end{equation}
\noindent where $\rho_{AB}$ is the Pearson correlation coefficient of $A$ and $B$.

However, this model can return misleading results when certain combinations of probabilities of events and correlations are employed. For example, considering $\mathbb{P}(A) = 0.3$ and $\mathbb{P}(B)=0.2$ with $\rho_{AB} = -1$ (opposite dependence), the Lucas model returns $\mathbb{P}(A \land B) = -0.123$, which is obviously erroneous.
The fact that Equation 2 is returning a negative probability is simply because the probability assigned for events A and B cannot have a correlation of -1.
This means that, for some probabilities of events, the Pearson correlation coefficient cannot take any value in [-1, 1], but in some subset $S \subseteq [-1,1]$.

This subset $S$ can be found through the Fréchet inequalities, which define the lower and upper bounds of the logical conjunction given the probability of its events \cite{frechet1935generalisation}, and are written as

\begin{equation}
    \mathbb{P}(A \land B) \in [\max(\mathbb{P}(A)+\mathbb{P}(B)-1, 0), \min(\mathbb{P}(A),\mathbb{P}(B))] \;.
\end{equation}

\noindent Substituting the bounds in Equation 2, and rearranging for $\rho_{AB}$, the subset $S = [\underline\rho_{AB},\overline\rho_{AB}] \subseteq [-1,1]$ can be found as

\begin{equation}
    \begin{cases}
      \underline{\rho}_{AB} = \frac{\max(\mathbb{P}(A)+\mathbb{P}(B)-1,0) - \mathbb{P}(A)\mathbb{P}(B)}{\sqrt{\mathbb{P}(A)\mathbb{P}(\neg A)\mathbb{P}(B)\mathbb{P}(\neg B)}}\\
      \\
      \overline{\rho}_{AB} = \frac{\min(\mathbb{P}(A),\mathbb{P}(B)) - \mathbb{P}(A)\mathbb{P}(B)}{\sqrt{\mathbb{P}(A)\mathbb{P}(\neg A)\mathbb{P}(B)\mathbb{P}(\neg B)}} \;.\\
    \end{cases}       
\end{equation}

\noindent With these definitions for the lower and upper bound of the Pearson correlation coefficient, the subset $S$ for the previous example with $\mathbb{P}(A) = 0.3$ and $\mathbb{P}(B) = 0.2$ is [-0.327, 0.763], and not [-1,1] as previously guessed.

The proposed correlated \textit{and} operation combines the Lucas model and the Fréchet inequalities to restrict the former to return probabilities in [0,1] for any specified Pearson correlation.
If the introduced correlation is greater than $\overline\rho$, then the upper Fréchet bound is returned.
On the other hand, if it is lower than $\underline\rho$, then the model gives the lower Fréchet bound.
Lastly, when the introduced correlation is in $S$, then the probability is calculated following the Lucas model.
Therefore, the correlated \textit{and} operator is written as

\begin{equation}
\mathbb{P}(A \land B) = \\
\begin{cases}
  \max(\mathbb{P}(A) + \mathbb{P}(B) -1, 0) \qquad \qquad \text{ if } \rho_{AB} \le \underline{\rho}_{AB} \\
  \min(\mathbb{P}(A),\mathbb{P}(B)) \qquad \qquad\qquad\qquad \text{ if } \rho_{AB} \ge \overline{\rho}_{AB} \\
  \mathbb{P}(A)\mathbb{P}(B) + \rho(A,B)\sqrt{\mathbb{P}(A)\mathbb{P}(\neg A)\mathbb{P}(B)\mathbb{P}(\neg B)} \text{ otherwise .}
\end{cases}
\end{equation}

We present an extension for the Boolean operators allowing for uncertainty not only in the inputs, but also in the dependence, which now can be specified for any interval correlation $\rho \subseteq [-1, 1]$, having the Fréchet bounds as a special case when $\rho = [-1,1]$. Functions of Bernoulli random variables with uncertainty in dependence generally yield interval probabilities, we thus show how intervals may also be propagated through the derived operations. 


\section{Correlated \textit{and} as a Copula}

\begin{figure}[!t]
    \centering
    \includegraphics[width=0.32\textwidth]{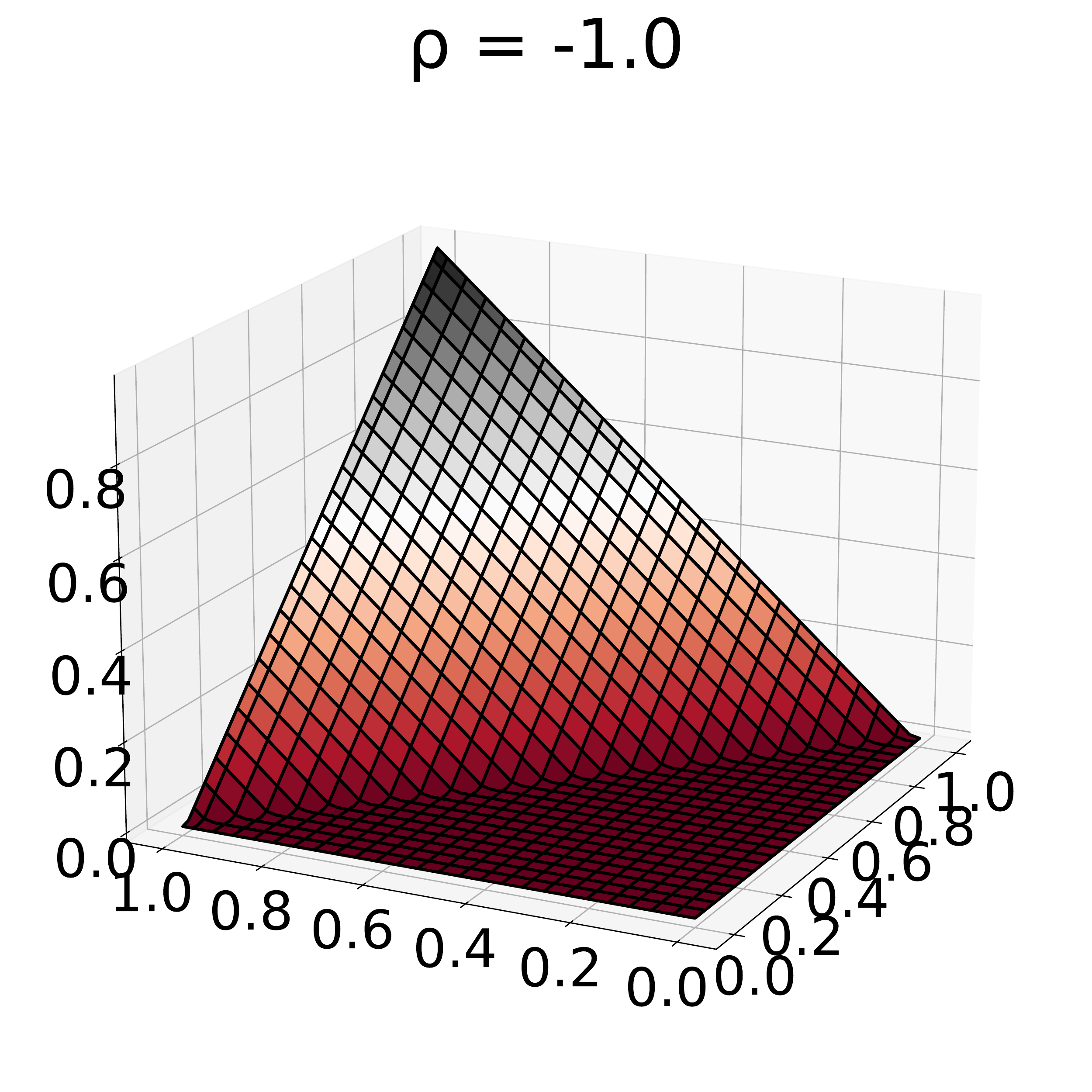}
    \includegraphics[width=0.32\textwidth]{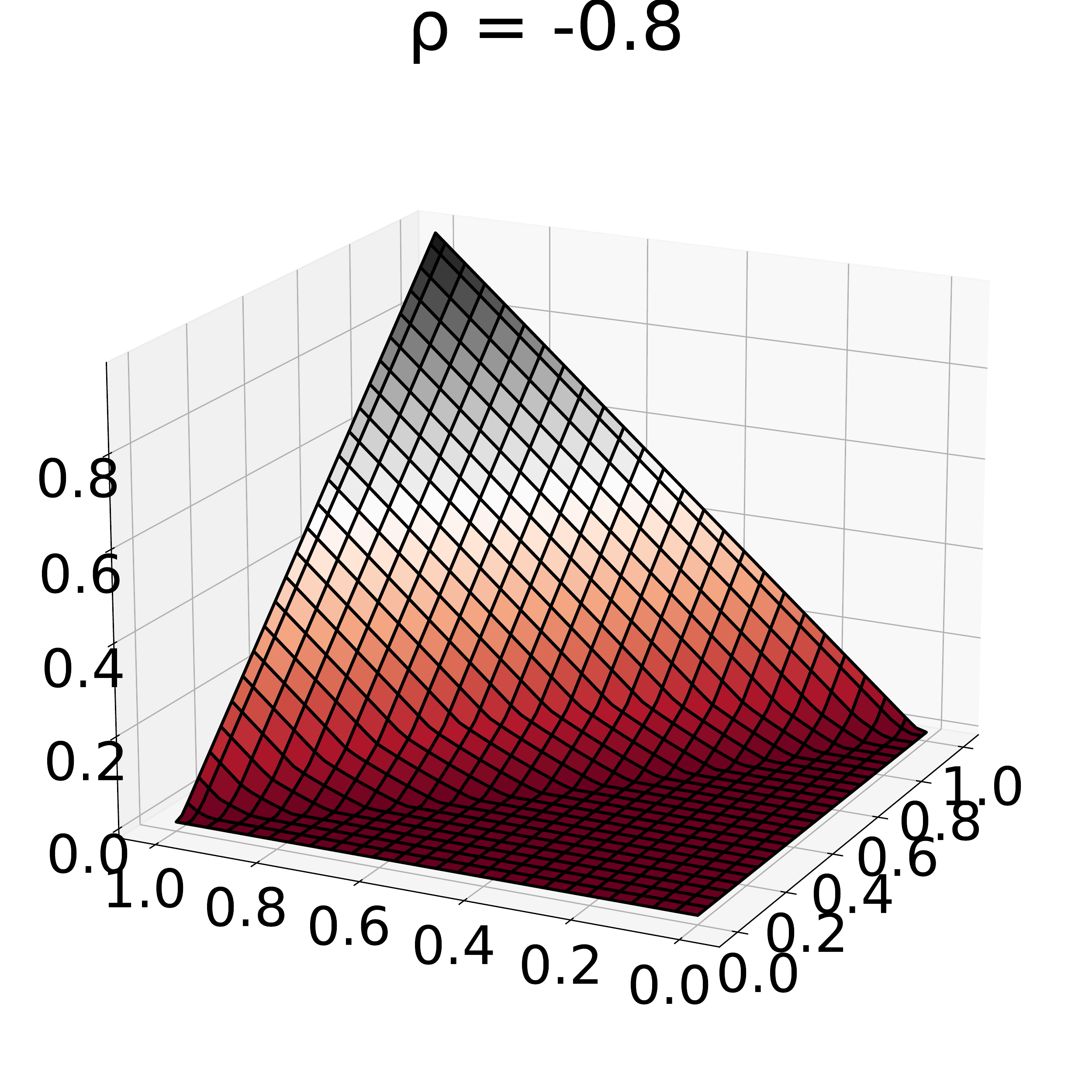}
    \includegraphics[width=0.32\textwidth]{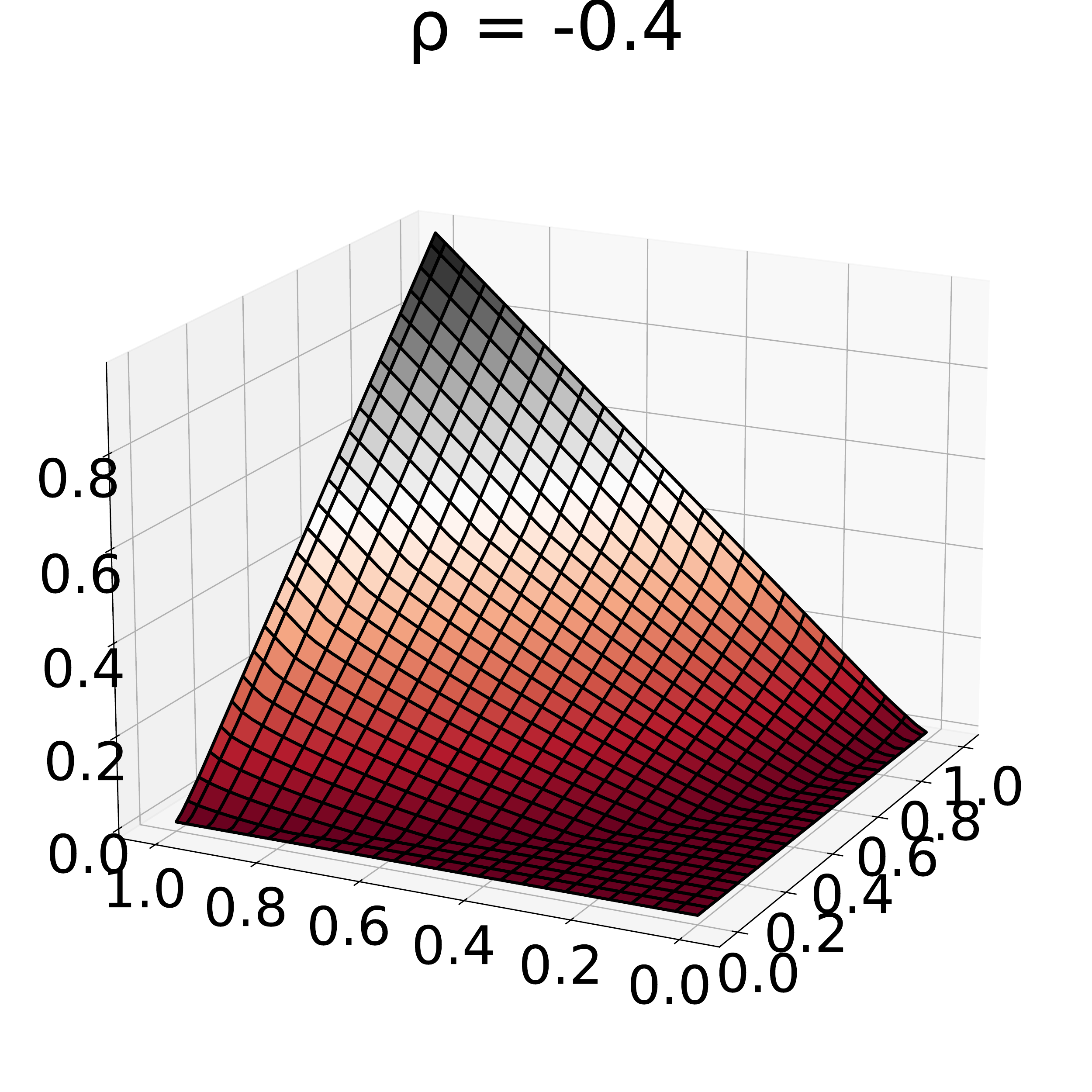}
    \includegraphics[width=0.32\textwidth]{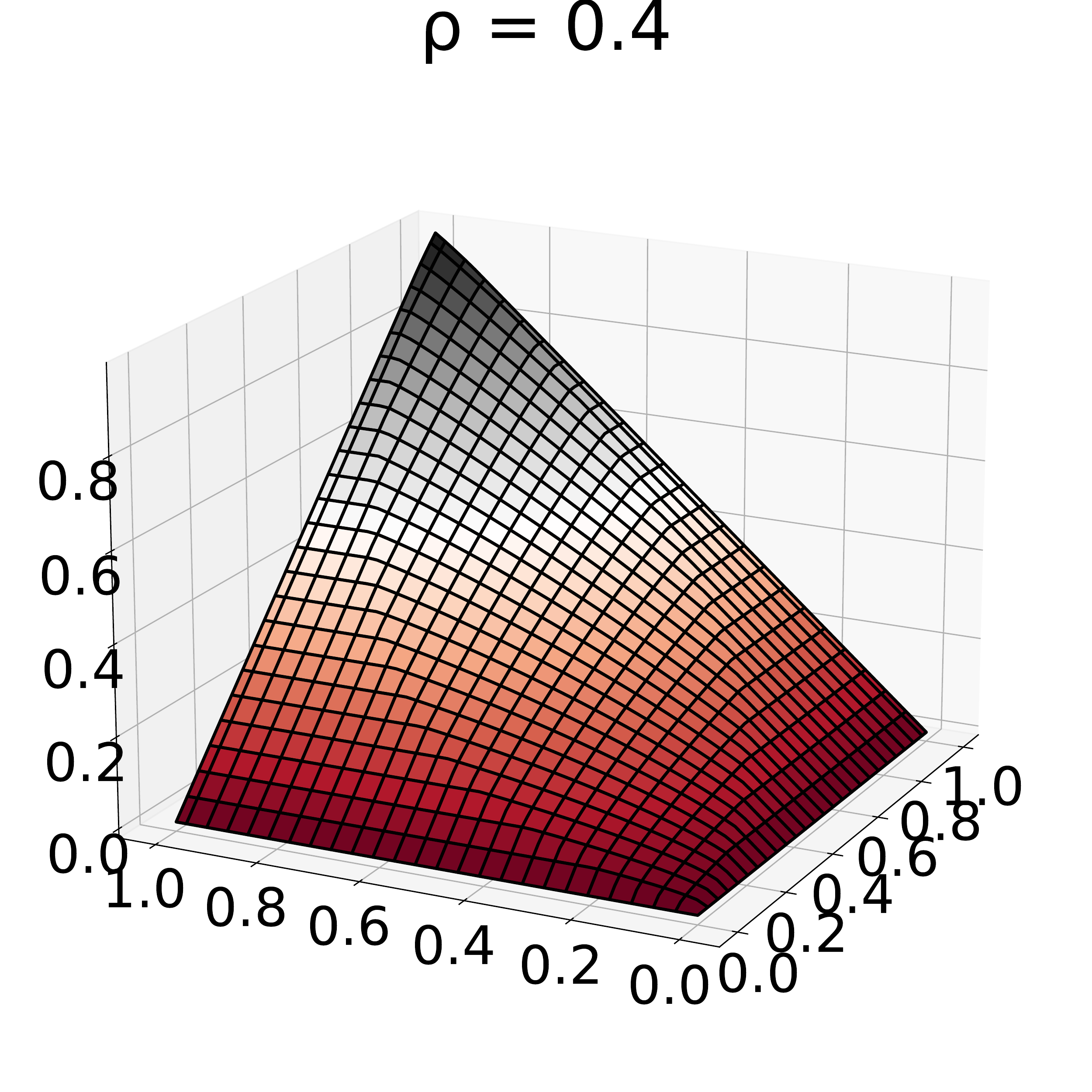}
    \includegraphics[width=0.32\textwidth]{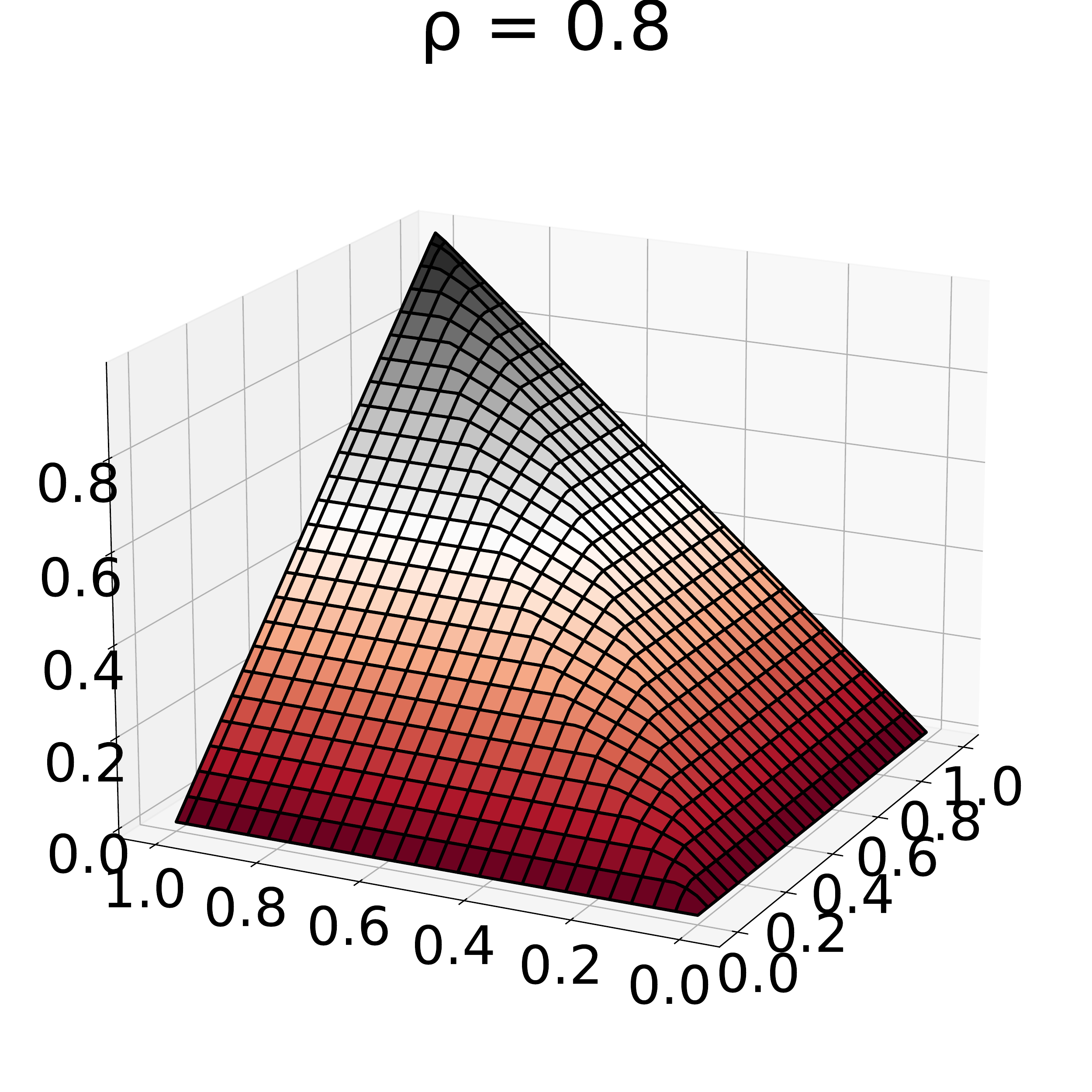}
    \includegraphics[width=0.32\textwidth]{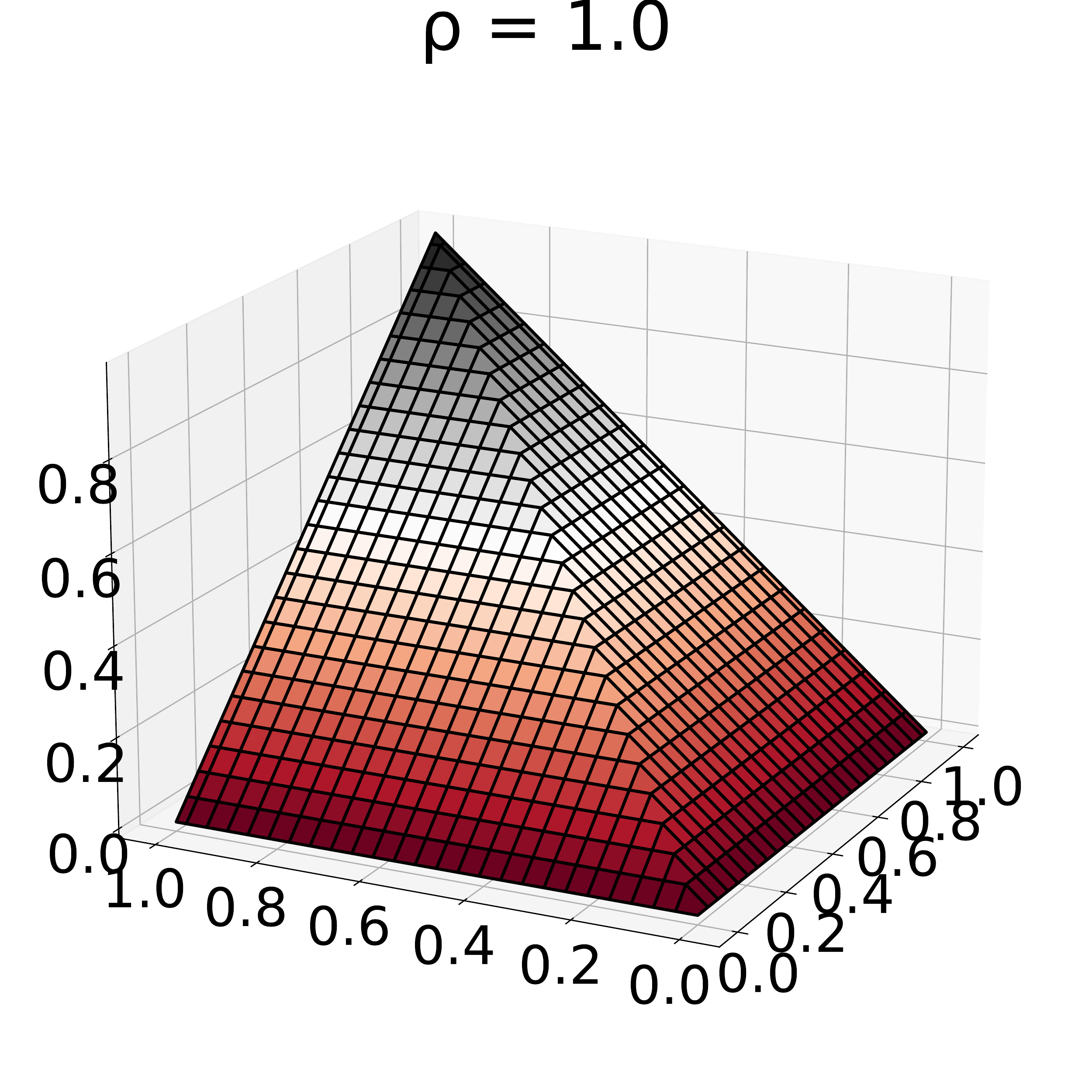}
    
    \caption{$C_\rho(u,v)$ for $\rho=\{-1, -0.8, -0.4, 0.4, 0.8, 1\}$}
    \label{fig:copulas}
\end{figure}

In this section we argue that the above derived correlated \textit{and} is a bivariate copula (2-copula) family, parameterised by a correlation coefficient $\rho$. Rewritten in a more standard copula notation:
\begin{equation*}
  C_{\rho}(u,v) =
    \begin{cases}
      W(u,v) & \text{if } \rho \le \underline{\rho}_{uv} \\
      M(u,v)& \text{if } \rho \ge \overline{\rho}_{uv} \\
      uv + \rho\sqrt{u(1- u)v(1- v)}& \text{otherwise} \;,
    \end{cases}
\end{equation*}

\noindent where $W(u,v) = \max(u + v - 1, 0)$ and $M(u,v) = \min(u, v)$ are the Fréchet-Hoeffding copula bounds. A 2-copula $C$ is any function $C:[0,1]^2 \rightarrow [0,1]$ with the following properties:

\begin{enumerate}
  \item Grounded: $C(0,v) = C(u,0) = 0$,
  \item Uniform margins: $C(u,1) = u; \;C(1,v) = v$,
  \item 2-increasing: \\ $C(u_{2},v_{2}) - C_{}(u_{2}, v_{1}) - C(u_{1}, v_{2}) + C(u_{1}, v_{1}) \ge 0$\\ for all $0 \le u_{1} \le u_{2} \le 1$ and $0 \le v_{1} \le v_{2} \le 1$.
\end{enumerate}
\noindent It is easy to see that the first two properties hold. The third property is harder to demonstrate, and in this work we provide no proof, despite the Lucas model $uv + \rho\sqrt{u(1- u)v(1- v)}$ is non-decreasing in $u$ and $v$. However, non-decreasing is a necessary but not sufficient criteria for 2-increasing \cite{schweizer2011probabilistic}. Yet Durante and Jaworski \cite{durante2010new}  prove that $C$ is a copula iff it satisfies criteria 1. and 2. and if the partial derivatives are increasing (Corollary 2.4). That is, for every $u \in [0,1]$,
\begin{equation*}
    v \mapsto \frac{\delta C(u,v)}{\delta u}
\end{equation*}
\noindent is increasing on $[0,1]$. The partial derivatives of the above \textit{and} gate is

\begin{equation*}
\frac{\delta C_{\rho}(u,v)}{\delta u} = 
\begin{cases}
  0 & \text{if } \rho_{uv} \le \underline{\rho}_{uv} \\
  1 & \text{if } \rho_{uv} \ge \overline{\rho}_{uv} \\
  v + \rho \frac{v(1-u)(1-v) - uv(1-v)}{2\sqrt{uv(1-u)(1-v)}} & \text{otherwise}
\end{cases}
\end{equation*}

\noindent which from observation $C_{\rho}$ follows. Figure \ref{fig:copulas} shows $C_{\rho}(u,v)$ for various values of $\rho$. If $C_{\rho}(u,v)$ is a copula family, then it is a \textit{complete} copula family, as it includes the two Fréchet-Hoeffding bounds $W$ and $M$, corresponding to minimal (when $\rho = -1$) and maximal (when $\rho = 1$) correlation respectively, and the independence copula $\Pi(u,v) = uv$ when $\rho = 0$.

\noindent For the main results of this paper it is not required that the presented operator is a copula, only that it is non-decreasing. However it is interesting since t-norms, which are similar functions to copulas, are widely used in fuzzy logic to model \textit{and} operations \cite{klir1987fuzzy}. The interpretation of t-norms in fuzzy logic is often unclear, whilst the presented $C_{\rho}$ has a clear probabilistic interpretation. Given two marginal Bernoulli random variables with $\mathbb{P}(A)$ and $\mathbb{P}(B)$, and Pearson correlation coefficient $\rho$, $C_{ \rho}(\mathbb{P}(A), \mathbb{P}(B))$ returns the joint probability of events $A$ and $B$ occurring, i.e. $\mathbb{P}(A=1, B=1)$, which is one element of the joint probability of the bivariate Bernoulli distribution, with the other elements being $\mathbb{P}(A=1, B=0)$, $\mathbb{P}(A=0, B=1)$, and $\mathbb{P}(A=0, B=0)$.

\section{Interval probabilities}
\label{sec:Interval_probs}
Since $C_{\rho}(u,v)$ is non-decreasing in $u$ and $v$, interval values for $u$ and $v$ can be simply evaluated with endpoints. Moreover, in $C_{\rho}(u,v)$ the correlation $\rho$ only has a single occurrence, and so it can be evaluated exactly with interval arithmetic. However, a useful observation is that $C_{\rho}(u,v)$ is also non-decreasing in $\rho$, and an interval value for $\rho$ induces an imprecise copula. The concept of imprecise copulas, which are a bounded set of copulas, has been discussed in \cite{montes2015sklar}. From Figure \ref{fig:copulas} it can be seen that increasing $\rho$ yields larger or equal values of joint probabilities, i.e. $C_{\rho_{1}} \leq C_{\rho_{2}}$ for $\rho_1 \leq \rho_2$, and therefore interval uncertainty in $\rho$ may also be evaluated with endpoints. Therefore, given two events with interval probabilities $u = [\underline{u}, \overline{u}]$ and $v =[\underline{v}, \overline{v}]$, and partially known correlation $\rho = [\underline{\rho}, \overline{\rho}]$, the interval probability of their conjunction can be evaluated as
\begin{align*}
    \underline{\mathbb{P}}(u=1, v = 1) &= C_{\underline{\rho}}(\underline{u}, \underline{v}) ,\\
    \overline{\mathbb{P}}(u=1, v = 1) &= C_{\overline{\rho}}(\overline{u}, \overline{v}) .
\end{align*}
\noindent As an example, for $u = [0.2,0.3]$, $v = [0.45,0.5]$ and $\rho =[ -0.2, 0.4]$, rigorous bounds on their conjunction can be calculated as $\mathbb{P}(u = 1, v =1) = [0.0502, 0.2417]$.

The full joint probability distribution of the bivariate Bernoulli can also be found by noticing that, for example, $\mathbb{P}(u =1, v = 0) = \mathbb{P}(u \land \neg v)$, with the \textit{not} operator defined as $\neg v = 1 - v$, which may be evaluated with $C_{\rho}$. However, some careful consideration is required regarding the correlation coefficient. If probabilities $u$ and $v$ have correlation $\rho_{uv}$, then $u$ and $1-v$ will have $\rho_{u\neg v} = -1 * \rho_{uv}$ \cite{nelsen2007introduction}, i.e. when complementing an event, the correlation must be reversed. The other elements of the joint distribution of the bivariate Bernoulli can be calculated as
\begin{align*}
   \mathbb{P}(u=1, v = 0) &= C_{-1*\rho}(u, 1 - v) ,\\
    \mathbb{P}(u=0, v = 1) &= C_{-1*\rho}(1 - u, v) ,\\
    \mathbb{P}(u=0, v = 0) &= C_{\rho}(1- u, 1 - v) .\\
\end{align*}
\noindent Notice that when the probabilities are complemented twice, the correlation stays the same, since it has been negated twice: $\rho_{\neg u, \neg v} = -(-\rho_{uv}) = \rho_{uv}$. Table \ref{tab:joint_bernoulli} shows the computed joint distribution using the previous example of $u = [0.2,0.3]$, $v = [0.45,0.5]$ and $\rho =[ -0.2, 0.4]$.

\begin{table}[t!]
\centering
\caption{Interval bounds on the joint probability of a bivariate Bernoulli with marginals $\mathbb{P}(u=1) = [0.2, 0.3]$ , $\mathbb{P}(v=1) = [0.45, 0.5]$ and interval correlation $\rho = [ -0.2, 0.4]$.}
\resizebox{0.6\textwidth}{!}{%
\begin{tabular}{c|cc|c}
\hline
$\mathbb{P}(u, v)$ & 0        & 1        &    $\mathbb{P}(u)$      \\ \hline
0           & $[0.304, 0.52]$ & $[0.223, 0.44]$ & $[0.7, 0.8]$ \\
1           & $[0.02, 0.2106]$ & $[0.0502, 0.2417]$ & $[0.2, 0.3]$ \\ \hline
$\mathbb{P}(v)$  & $[0.5, 0.55]$ & $[0.45, 0.5]$ &        \\ \hline
\end{tabular}%
}
\label{tab:joint_bernoulli}
\end{table}

\subsection{Other Boolean operations}

Once a correlated $and$ operation has been constructed, other correlated Boolean operations can be defined in terms of this operation and a \textit{not}, shown in Table \ref{Tabel1}. For example, \textit{or} can be defined as

\begin{equation*}
    A \text{ or } B = \mathbb{P}(A \lor B) = 1 - \mathbb{P}((1 -A) \land (1-B)),
\end{equation*}

\noindent which may be written in terms of $C_{\rho}$ is

\begin{equation}
    A \lor B = 1 - C_{\rho}(1-A,1-B).
\end{equation}

\noindent The other operations can be similarly expanded, taking care to negate $\rho$ appropriately.

\section{Probability Boxes}

The previous section describes a method to perform logical operations on uncertain Booleans characterised by interval probabilities $p \subseteq [0, 1]$, which define a bounded set of Bernoulli distributions. A possible generalisation of this is to have a distributional or imprecise distributional (p-box) characterisation, e.g., any p-box $p$ whose range is a subset of the unit interval $\text{range}(p)\subseteq [0,1]$. The binary events involved with Boolean operations can be expressed in the form of Bernoullis (precise), set of Bernoullis (interval), or distributional Bernoullis (p-box). Note that here we are not describing the arithmetic of real functions ($f:R^{m} \to R^{n}$), but are describing an extension of Boolean functions with p-box inputs, i.e., events or uncertain Booleans with uncertainty characterised by p-boxes. This is particularly relevant for c-boxes \cite{ferson2014computing}, which are p-box shaped confidence distributions for binomial inference with limited data. That is, they are a confidence characterisation of an uncertain Boolean given some sample set, e.g., the data for $A$ in Equation \ref{eq:rand_vec}.

\begin{figure}[t!]
    \centering
\tikzset{every picture/.style={line width=0.75pt}} 
\begin{tikzpicture}[x=0.75pt,y=0.75pt,yscale=-1,xscale=1]
\draw (216,185) node  {\includegraphics[width=180pt,height=172.5pt]{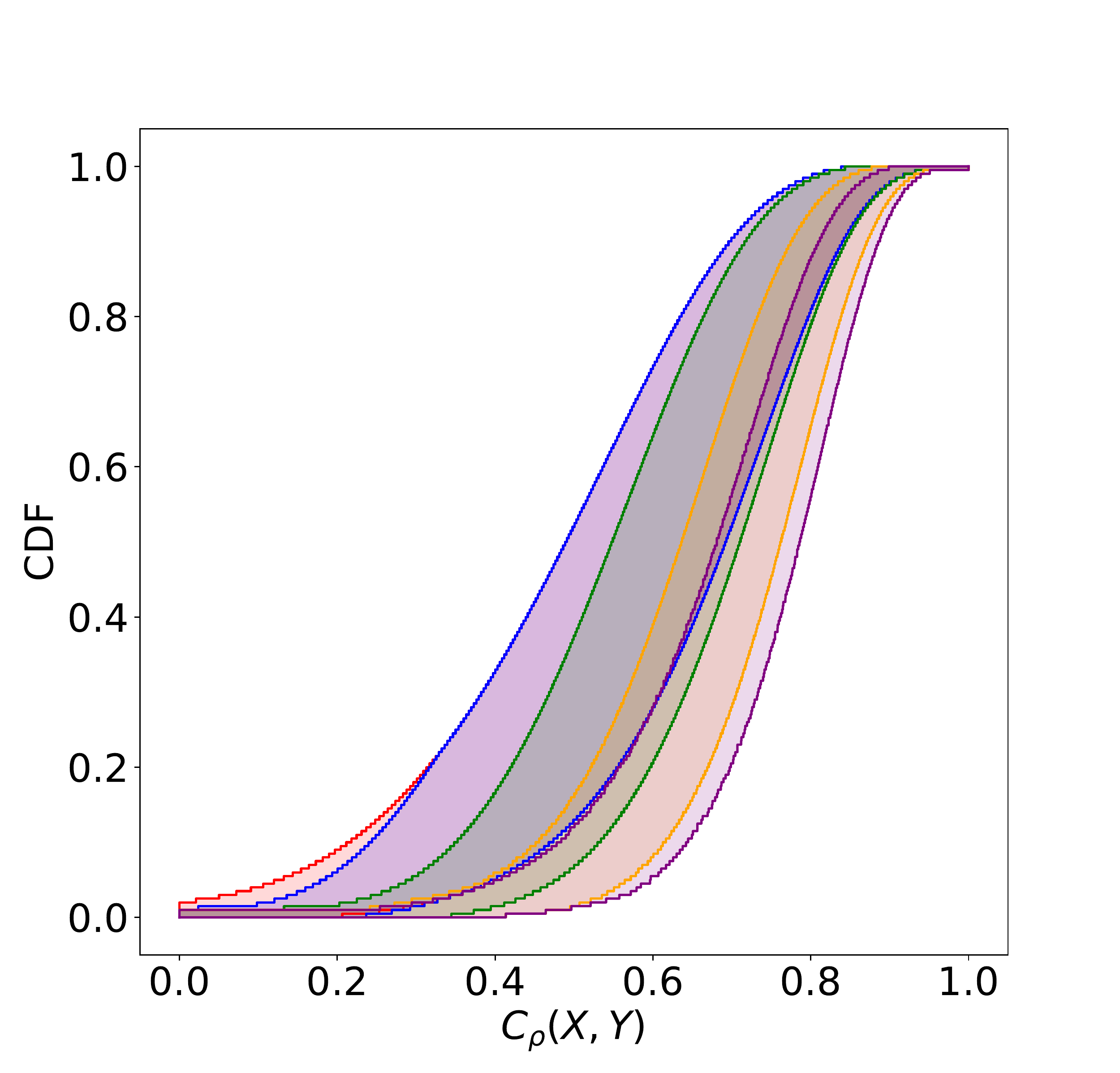}};
\draw (450,185) node  {\includegraphics[width=180pt,height=172.5pt]{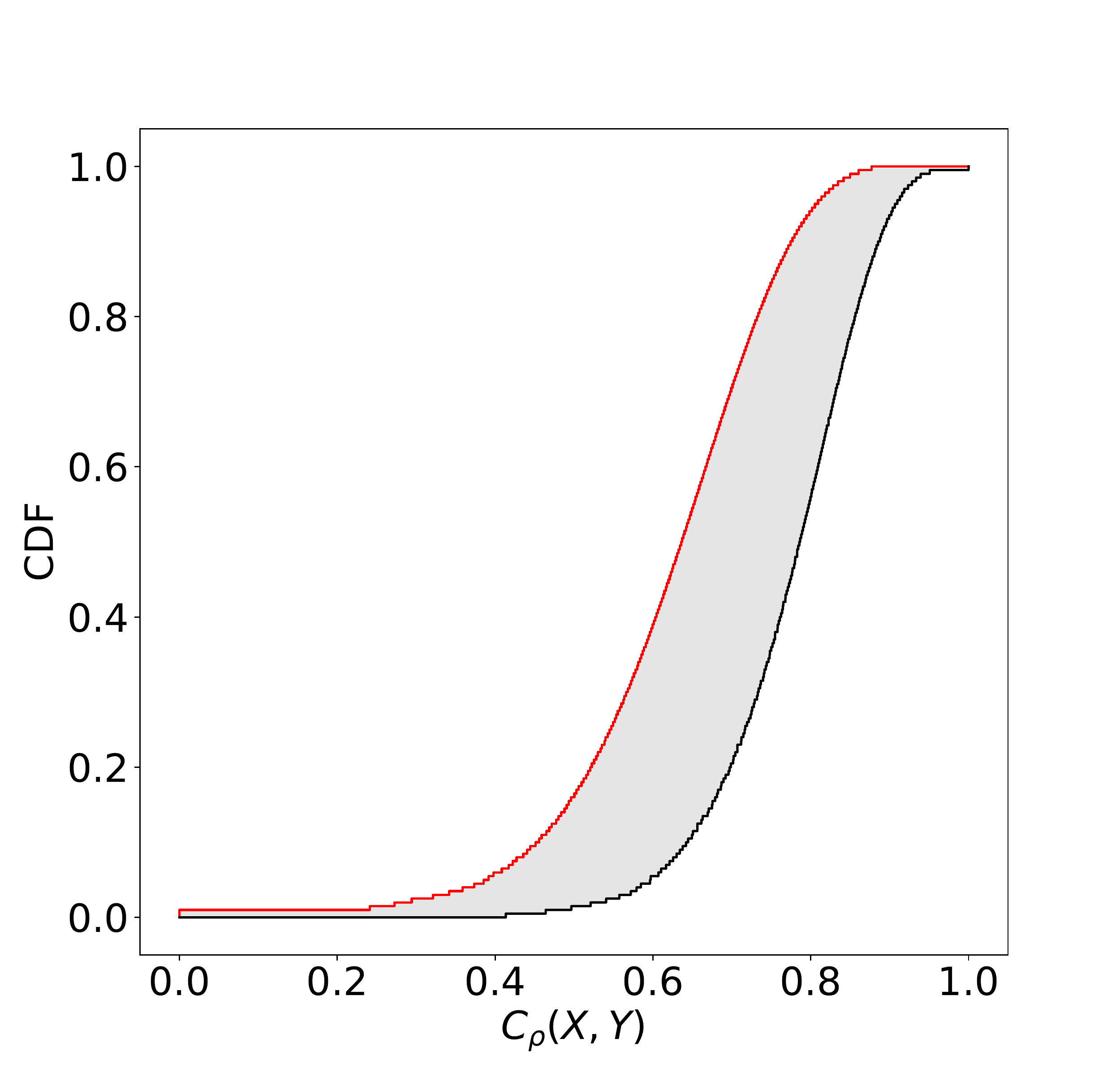}};
\draw [color={rgb, 255:red, 255; green, 0; blue, 30 }  ,draw opacity=1 ][line width=1.5]    (156,90) .. controls (151.52,115.8) and (138.15,137.3) .. (136,150) .. controls (133.88,162.51) and (115.89,201.31) .. (154.21,257.42) ;
\draw [shift={(156,260)}, rotate = 234.73] [fill={rgb, 255:red, 255; green, 0; blue, 30 }  ,fill opacity=1 ][line width=0.08]  [draw opacity=0] (13.4,-6.43) -- (0,0) -- (13.4,6.44) -- (8.9,0) -- cycle    ;
\draw [color={rgb, 255:red, 39; green, 110; blue, 193 }  ,draw opacity=1 ][line width=1.5]    (176,120) .. controls (152.69,159.89) and (149.43,168.25) .. (184.38,246.38) ;
\draw [shift={(186,250)}, rotate = 245.78] [fill={rgb, 255:red, 39; green, 110; blue, 193 }  ,fill opacity=1 ][line width=0.08]  [draw opacity=0] (13.4,-6.43) -- (0,0) -- (13.4,6.44) -- (8.9,0) -- cycle    ;
\draw [color={rgb, 255:red, 78; green, 133; blue, 17 }  ,draw opacity=1 ][line width=1.5]    (196,150) .. controls (193.41,180.56) and (192.39,190.88) .. (213.94,226.61) ;
\draw [shift={(216,230)}, rotate = 238.42] [fill={rgb, 255:red, 78; green, 133; blue, 17 }  ,fill opacity=1 ][line width=0.08]  [draw opacity=0] (13.4,-6.43) -- (0,0) -- (13.4,6.44) -- (8.9,0) -- cycle    ;
\draw [color={rgb, 255:red, 197; green, 124; blue, 1 }  ,draw opacity=1 ][line width=1.5]    (296,200) .. controls (295.35,173.33) and (305.68,157.47) .. (269.49,168.86) ;
\draw [shift={(266,170)}, rotate = 341.49] [fill={rgb, 255:red, 197; green, 124; blue, 1 }  ,fill opacity=1 ][line width=0.08]  [draw opacity=0] (13.4,-6.43) -- (0,0) -- (13.4,6.44) -- (8.9,0) -- cycle    ;
\draw [color={rgb, 255:red, 144; green, 19; blue, 254 }  ,draw opacity=1 ][line width=1.5]    (296,250) .. controls (300.16,230.32) and (287.42,235.06) .. (259.56,230.61) ;
\draw [shift={(256,230)}, rotate = 10.28] [fill={rgb, 255:red, 144; green, 19; blue, 254 }  ,fill opacity=1 ][line width=0.08]  [draw opacity=0] (13.4,-6.43) -- (0,0) -- (13.4,6.44) -- (8.9,0) -- cycle    ;
\draw (134,72.4) node [anchor=north west][inner sep=0.75pt]    {$\textcolor[rgb]{1,0,0}{\rho \ =\ -1}$};
\draw (157,102.4) node [anchor=north west][inner sep=0.75pt]    {$\textcolor[rgb]{0.15,0.43,0.76}{\rho \ =\ -0.5}$};
\draw (175,132.4) node [anchor=north west][inner sep=0.75pt]    {$\textcolor[rgb]{0.31,0.54,0.07}{\rho \ =\ 0}$};
\draw (264,202.4) node [anchor=north west][inner sep=0.75pt]    {$\textcolor[rgb]{0.77,0.49,0}{\rho \ =\ 0.5}$};
\draw (272,252.4) node [anchor=north west][inner sep=0.75pt]    {$\textcolor[rgb]{0.56,0.07,1}{\rho \ =\ 1}$};
\draw (400,72.4) node [anchor=north west][inner sep=0.75pt]  [font=\large]  {$\textcolor[rgb]{0,0,0}{\rho \ =\ [ 0.5,\ 1]}$};
\end{tikzpicture}
    \caption{\textbf{(Left)} Show logical \textit{and} between the two K-out-of-N c-boxes $X \sim \text{KN}(5,6)$, $Y \sim \text{KN}(16,20)$ for different values of precise correlation. \textbf{(Right)} shows logical \textit{and} with the interval correlation $\rho = [0.5,1]$, and is the envelope of the orange and purple c-boxes. Independence is used in the upper level throughout $C_{XY} = \Pi$.}
    \label{fig:pbox_ands}
\end{figure}
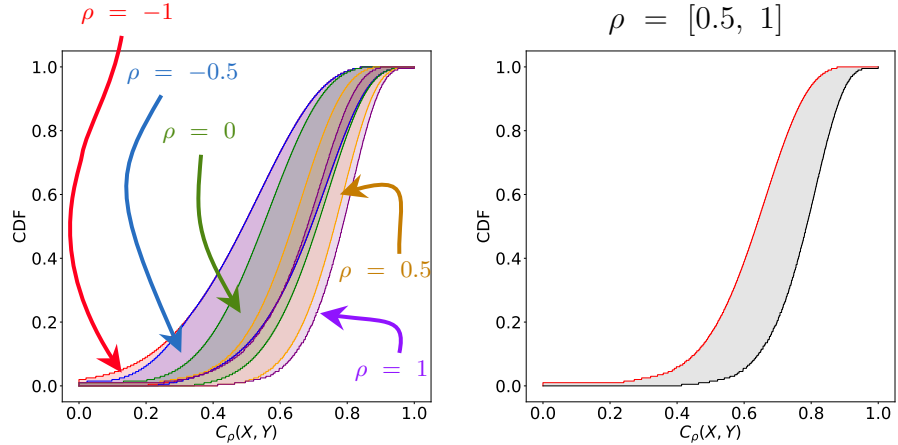
Since $C_{\rho}$ is a non-decreasing binary operator, it can also be readily evaluated with the convolutions used in p-box arithmetic. Given two random variables with distribution functions $F_{X}$ and $F_{Y}$, correlated by copula $C_{XY}$, a binary operation $Z= f(X,Y)$ can be evaluated with the following Lebesgue-Stieltjes integral 
\begin{equation*}
  F_{Z}(z) = \int_{f\{z\}}dC_{XY}(F_{X}(x),F_{Y}(y)) ,
\end{equation*}
\noindent where the integration domain is the set $f\{z\} = \{(x,y)| x,y \in \mathbb{R}, f(x,y) < z\}$ (the set of all $x$ and $y$ for which $f(x,y) < z$). Note that the copula $C_{XY}$ is not the same copula used to define the correlated \textit{and} operation $C_{\rho}$, where $C_{XY}$ defines the dependence between the random variables $X$ and $Y$ (which defines the uncertainty we have about the events), and $C_{\rho}$ defines the correlation between the events themselves. We elaborate on this difference in Section \ref{sec:2_level}.

\noindent Inserting $f = C_{\rho}$ into the above equation, and since the operation is non-decreasing we have
\begin{align}
    \underline{F_{Z}}(z) &= \int_{C_{\underline{\rho}}}dC_{XY}\left(\underline{F}_{X}(x), \underline{F}_{Y}(y)\right)
    , \label{eq:pbox_and1} \\
    \overline{F_{Z}}(z) &= \int_{C_{\overline{\rho}}}dC_{XY}\left(\overline{F}_{X}(x), \overline{F}_{Y}(y)\right) .
    \label{eq:pbox_and2}
\end{align}
\noindent for two p-boxes $X = [\underline{F}_{X}, \overline{F}_{X}]$, $Y = [\underline{F}_{Y}, \overline{F}_{Y}]$, and interval correlation $\rho = [\underline{\rho}, \overline{\rho}]$. Software for performing rigorous correlated p-box arithmetic is readily available for bounding the above integrals efficiently \cite{gray2021probabilityboundsanalysis}. The above integrals are those usually used in p-box arithmetic, except that the binary operation $C_{\rho}$ is parameterised by $\rho$, which could be an interval. In the interval case, the envelope of the two end points yields the bounds on the output p-box. Figure \ref{fig:pbox_ands} shows Equations \ref{eq:pbox_and1} and  \ref{eq:pbox_and2} evaluated for two K-out-of-N c-boxes \cite{ferson2014computing} $X \sim \text{KN}(5,6)$, $Y \sim \text{KN}(16,20)$ for different values of event correlations, shown in different colours. The right figure shows logical \textit{and} with an interval correlation $\rho = [0.5,1]$, which is the envelope of the purple and orange c-boxes. Note in Figure \ref{fig:pbox_ands} independence was used for the upper level $C_{XY}$, but independence is not necessarily the only choice for $C_{XY}$. In fact, any copula $C_{XY}$ can be used in Equations \eqref{eq:pbox_and1} and \eqref{eq:pbox_and2}, which can greatly influence results, as much as the choice of correlation $\rho$ does. We explore this difference in the next section.

\subsection{Two levels of dependence}
\label{sec:2_level}
The copula $C_{XY}$ in the \textit{and} operation \eqref{eq:pbox_and1} and \eqref{eq:pbox_and2} plays a different role to the correlation $\rho$. The correlation $\rho$ is the dependence between the events $\mathbb{P}(A) = a$ and $\mathbb{P}(B) = b$, i.e., the correlation between the two random bit-vectors \eqref{eq:rand_vec}, and is the only dependence that plays a role when two uncertain Booleans are characterised by precise probabilities as real values. When considering two interval probabilities, two distributions, or two p-boxes, $\rho$ still correlates the events as before, however one could worry that the uncertainty characterising the marginal probabilities might share some bivariate information. That is, the p-boxes have a dependence which the copula $C_{XY}$ characterises, distinctly from the event correlation $\rho$. 
\begin{figure}[t!]
    \centering
    \tikzset{every picture/.style={line width=0.75pt}} 
    
    \begin{tikzpicture}[x=0.75pt,y=0.75pt,yscale=-1,xscale=1]
    
    \draw (216,185) node  {\includegraphics[width=180pt,height=172.5pt]{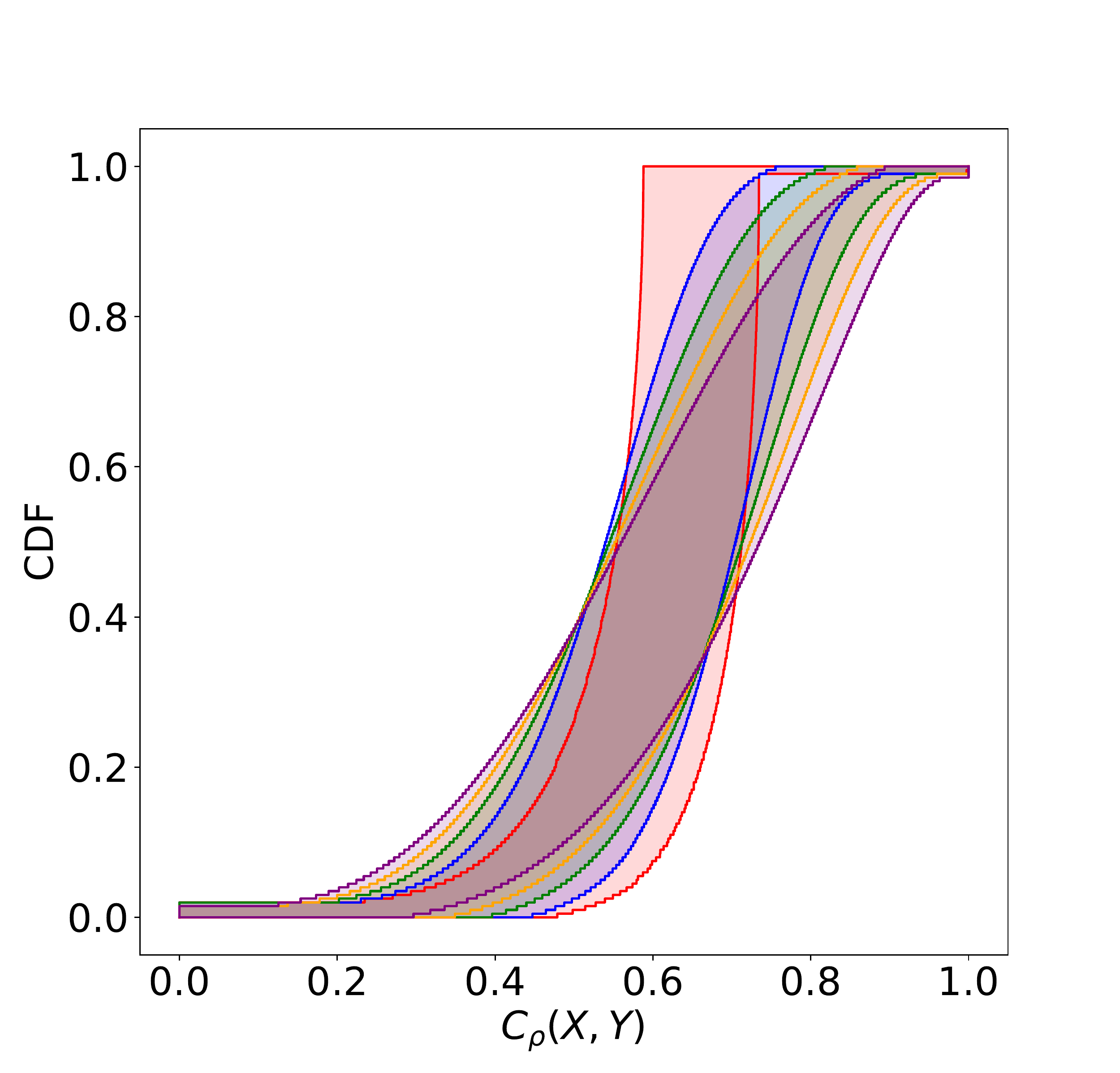}};
    \draw (450,185) node  {\includegraphics[width=180pt,height=172.5pt]{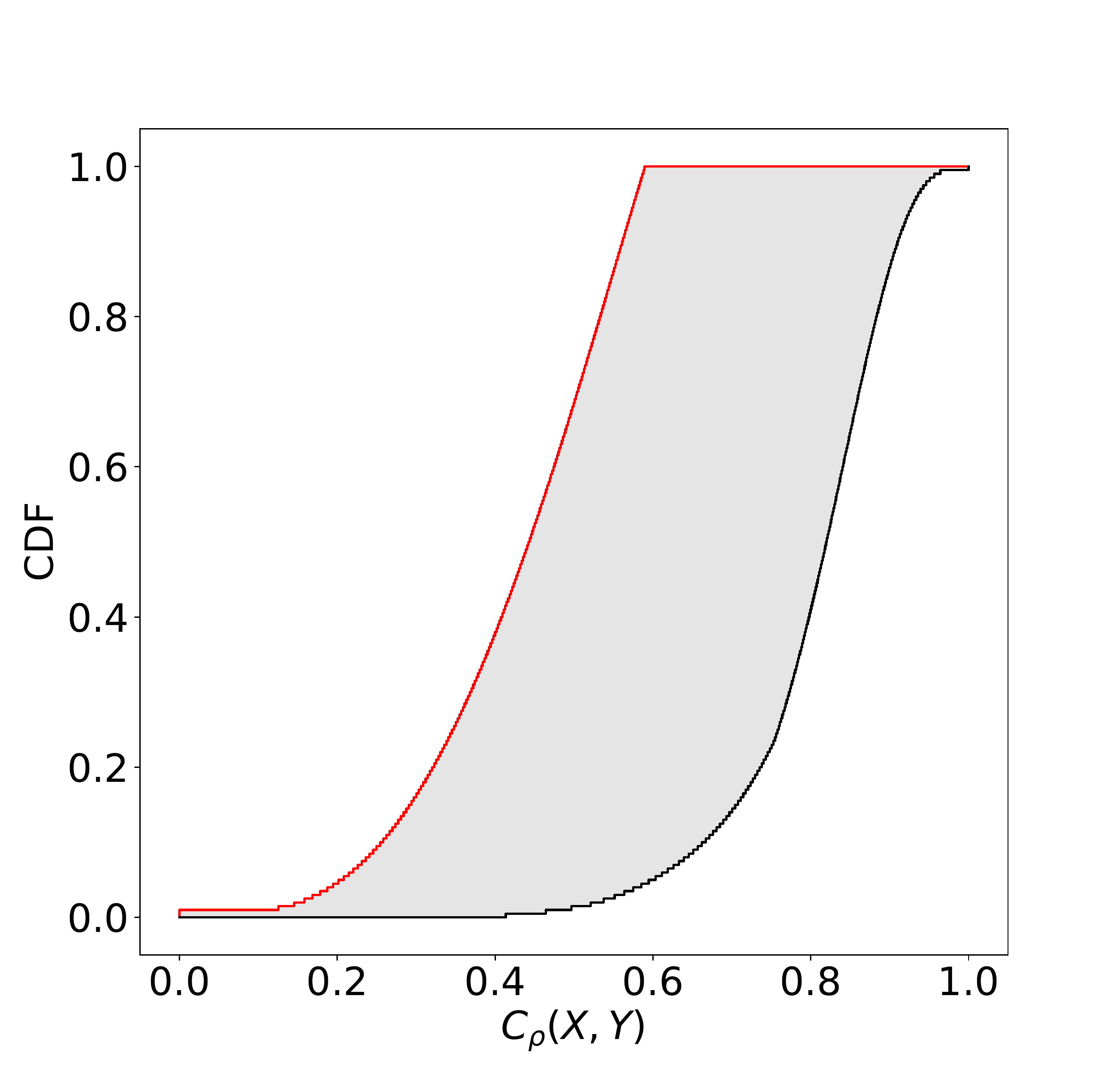}};
    \draw [color={rgb, 255:red, 255; green, 0; blue, 30 }  ,draw opacity=1 ][line width=1.5]    (170,110) .. controls (184.08,121.59) and (178.1,134.56) .. (226.19,130.36) ;
    \draw [shift={(230,130)}, rotate = 174.29] [fill={rgb, 255:red, 255; green, 0; blue, 30 }  ,fill opacity=1 ][line width=0.08]  [draw opacity=0] (13.4,-6.43) -- (0,0) -- (13.4,6.44) -- (8.9,0) -- cycle    ;
    \draw [color={rgb, 255:red, 39; green, 110; blue, 193 }  ,draw opacity=1 ][line width=1.5]    (160,160) .. controls (180.37,179.84) and (195.98,155.7) .. (226.61,168.46) ;
    \draw [shift={(230,170)}, rotate = 206.1] [fill={rgb, 255:red, 39; green, 110; blue, 193 }  ,fill opacity=1 ][line width=0.08]  [draw opacity=0] (13.4,-6.43) -- (0,0) -- (13.4,6.44) -- (8.9,0) -- cycle    ;
    \draw [color={rgb, 255:red, 78; green, 133; blue, 17 }  ,draw opacity=1 ][line width=1.5]    (170,200) .. controls (177.58,218.09) and (154.32,226.15) .. (196.62,239) ;
    \draw [shift={(200,240)}, rotate = 196] [fill={rgb, 255:red, 78; green, 133; blue, 17 }  ,fill opacity=1 ][line width=0.08]  [draw opacity=0] (13.4,-6.43) -- (0,0) -- (13.4,6.44) -- (8.9,0) -- cycle    ;
    \draw [color={rgb, 255:red, 197; green, 124; blue, 1 }  ,draw opacity=1 ][line width=1.5]    (290,220) .. controls (272.23,184.79) and (266.72,200.39) .. (269.75,153.72) ;
    \draw [shift={(270,150)}, rotate = 93.97] [fill={rgb, 255:red, 197; green, 124; blue, 1 }  ,fill opacity=1 ][line width=0.08]  [draw opacity=0] (13.4,-6.43) -- (0,0) -- (13.4,6.44) -- (8.9,0) -- cycle    ;
    \draw [color={rgb, 255:red, 144; green, 19; blue, 254 }  ,draw opacity=1 ][line width=1.5]    (320,170) .. controls (324.16,150.32) and (317.06,137.6) .. (293.08,121.97) ;
    \draw [shift={(290,120)}, rotate = 32.05] [fill={rgb, 255:red, 144; green, 19; blue, 254 }  ,fill opacity=1 ][line width=0.08]  [draw opacity=0] (13.4,-6.43) -- (0,0) -- (13.4,6.44) -- (8.9,0) -- cycle    ;
    
    \draw (131,98.4) node [anchor=north west][inner sep=0.75pt]    {$\textcolor[rgb]{1,0,0}{C_{XY} =W}$};
    \draw (127,141.4) node [anchor=north west][inner sep=0.75pt]    {$\textcolor[rgb]{0.15,0.43,0.76}{C_{X}{}_{Y} =\Phi _{-0.5}}$};
    \draw (132,181.4) node [anchor=north west][inner sep=0.75pt]    {$\textcolor[rgb]{0.31,0.54,0.07}{C_{XY} =\Pi }$};
    \draw (256,221.4) node [anchor=north west][inner sep=0.75pt]    {$\textcolor[rgb]{0.77,0.49,0}{C_{XY} =\Phi _{0.5}}$};
    \draw (277,171.4) node [anchor=north west][inner sep=0.75pt]    {$\textcolor[rgb]{0.56,0.07,1}{C_{XY} =M}$};
    \draw (395,72.4) node [anchor=north west][inner sep=0.75pt]  [font=\large]  {$\textcolor[rgb]{0,0,0}{C_{XY} \ =\ }\textcolor[rgb]{0,0,0}{[}\textcolor[rgb]{0,0,0}{W,\ M}\textcolor[rgb]{0,0,0}{]}$};

    \end{tikzpicture}
    \caption{Shows the conjunction of two p-boxes with varying upper-level dependence. \textbf{(Left)} Show logical \textit{and} between the two K-out-of-N c-boxes $X \sim \text{KN}(5,6)$, $Y \sim \text{KN}(16,20)$ for constant independent event correlation $\rho$, but different values of upper-level dependence. $\Phi_r$ is the Gaussian copula with parameter value $r$. \textbf{(Right)} shows logical \textit{and} with independence for the lower-level $\rho$, but unknown dependence in the upper-level.}
    \label{fig:pbox_Frechet}
\end{figure}
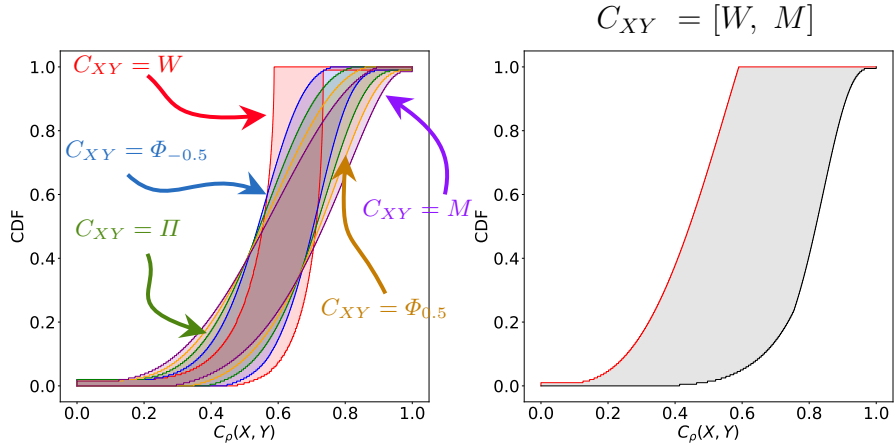
We call the event correlation $\rho$ as ``lower-level dependence'' and the copula 
$C_{XY}$ as the ``upper-level dependence''. This upper-level dependence can have a profound impact on the conjunction of two uncertain Booleans. The left of Figure \ref{fig:pbox_Frechet} shows the variation in the same two c-boxes as Figure \ref{fig:pbox_ands}, but with a constant lower-level independence and a varying upper-level $C_{XY}$. The copula $\Phi_{r}$ is a Gaussian copula with parameter $r$. Note that the green p-box, showing independence on both levels, is the same in both Figures \ref{fig:pbox_ands} and \ref{fig:pbox_Frechet}. 

As for the lower-level, the dependence on the upper-level could be imprecisely known, that is, it may not be possible to know $C_{XY}$. Using p-box arithmetic, unknown dependence can be propagated through a non-decreasing binary operation using the following convolutions \cite{ferson2015dependence} which has been adapted using the imprecise probabilistic conjunction $C_{\rho}$
\begin{align}
    \underline{F_{Z}}(z) &= \inf_{C_{\underline{\rho}}(x,y) = z} \left[W\left(\underline{F}_{X}(x), \underline{F}_{Y}(y)\right) \right]
    , \label{eq:pbox_tau1} \\
    \overline{F_{Z}}(z) &= \sup_{C_{\overline{\rho}}(x,y) = z}\left[W^{d}\left(\overline{F}_{X}(x), \overline{F}_{Y}(y)\right) \right],
    \label{eq:pbox_tau2}
\end{align}
\noindent where $C^{d}$ is the dual copula of $C$: $C^{d}(u,v) = u+v-C(u,v)$. The right of Figure \ref{fig:pbox_Frechet} shows the result of the p-box conjunction with unknown upper-level dependence and with lower-level independence. 

The p-box bounds from unknown upper-level dependence are not only rigorous but are also best-possible in the sense that they cannot be made tighter without introducing additional dependence assumptions. One may find the breath (grey shaded area) of this p-box surprising, and how strong the independence assumption is. One may also wonder why the p-box on the right is not a simple envelope of the p-boxes on the left. This is because the left p-boxes only consider Gaussian dependencies, whilst the right p-box considers all possible copulas $C_{XY}$, of which there are an infinite number.

Risk analysts tell us that for a fault tree analysis to be probabilistic the event probabilities should be characterised by distributions \cite{atwood2003handbook}.
If so, then the issue of the two levels of dependence immediately arises.
As far as we are aware, this issue has not been widely addressed.
Generally, the upper-level independence is assumed.
This section has shown that this assumption has a strong effect on the results of logical operations involving uncertainty.


\section{Application}

In fault tree analysis, Boolean operators are employed to calculate a system's probability of failure.
These are useful to, for example, understand what are the most likely routes of failure, show compliance with the reliability requirements, or designing monitoring strategies.
To carry out these analyses it is required to define an event tree, where the connections between events are made with Boolean operators, with their respective event probabilities and dependencies.
Then, the Boolean operators are calculated backwards from the top event to find its probability of failure.

Figure 4 represents a fault tree for a pressure tank system, derived in \cite{vesely1981fault,ferson2015dependence}.
It displays the Boolean operators $\land$ (AND) and $\lor$ (OR) connecting the failure events $E_i$ (with $E_1$ being the top event), and the system components tank (T), relay (K2), pressure switch (S), on-switch (S1), timer relay (R), and relay (K1).
For the sake of brevity, the function of these components will not be explained in this paper.

\begin{figure}[!ht]
    \centering
    \includegraphics[width=0.8\textwidth]{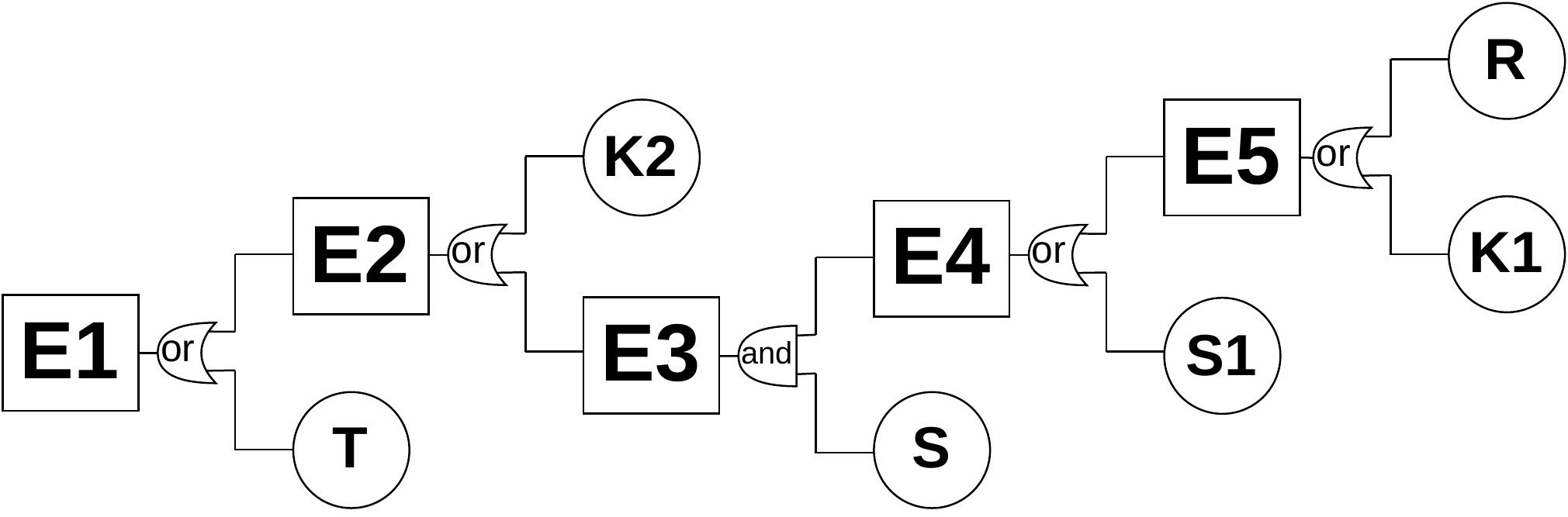}   \caption{Fault tree diagram for the pressure tank system.}
    \label{fig:my_label}
\end{figure}

To measure the probability of failure of  the system (event $E_1$), it is necessary to estimate the failure rates for its components T, K2, S, S1, R and K1, and the dependence of the events in the Boolean operations.
The most common method of calculation, as in \cite{vesely1981fault}, is to calculate the probability of failure assuming independence and known rates of failure for the components.
In \cite{ferson2015dependence}, the assumption of lower-level independence was relaxed including some known/unknown dependence, and extending the point probabilities to intervals, showing that the assumptions in \cite{vesely1981fault} underestimate the probability of failure.
To demonstrate the capabilities of the method derived in this paper, two scenarios on the components probability of failure will be calculated, with three lower-level dependence cases for each (independence, mixed dependence, and unknown dependence):

\begin{enumerate}
    \item An interval scenario, where the probability of failure of the components is modeled as intervals.
    \item A p-box scenario, which generalises \cite{ferson2015dependence} combining p-boxes and intervals. The probabilities of failure of the relays (K1, K2, and R) have been extended to follow a K-out-of-N c-box \cite{ferson2014computing}. 
\end{enumerate}

\noindent The specific values for each scenario and dependence are indicated in Tables 3 and 4 respectively.

\begin{table}[]
\centering
\caption{Probability of failure for the tank pressure system components. In the interval scenario, all probabilities are modeled in the form of intervals. In the p-box scenario, tank and switches have an interval probability, whilst the relays follow a K-out-of-N c-box.}\label{tab1}
\begin{tabular}{|c|c|c|}
\hline
Component & Interval scenario            & P-box scenario                           \\ \hline
T         & ${[}4.5\times10^{-6}, 5.5\times10^{-6}{]}$ & ${[}4.5\times10^{-6}, 5.5\times10^{-6}{]}$              \\ \hline
K2        & ${[}2.5\times10^{-5}, 3.5\times10^{-5}{]}$ & $KN(3,10^5)$ \\ \hline
S         & ${[}0.5\times10^{-4}, 1.5\times10^{-4}{]}$ & ${[}0.5\times10^{-4}, 1.5\times10^{-4}{]}$              \\ \hline
K1        & ${[}2.5\times10^{-5}, 3.5\times10^{-5}{]}$ & $KN(3,10^5)$ \\ \hline
R         & ${[}0.5\times10^{-4}, 1.5\times10^{-4}{]}$ & $KN(1,10^4)$ \\ \hline
S1        & ${[}2.5\times10^{-5}, 3.5\times10^{-5}{]}$ & ${[}2.5\times10^{-5}, 3.5\times10^{-5}{]}$              \\ \hline
\end{tabular}
\end{table}

\begin{table}[!ht]
\centering
\caption{Dependence in the system's events for the three different cases: independence, mixed dependence, and unknown dependence.}\label{tab1}
\begin{tabular}{|c|c|c|c|}
\hline
 Event  & Independence & Mixed Dependence                         & Unknown Dependence                       \\ \hline
$E_1$ & $\rho = 0$   & $\rho = 0$                               & $\rho = [-1, 1]$ \\ \hline
$E_2$ & $\rho = 0$   & $\rho = [-1, 1]$ & $\rho = [-1, 1]$ \\ \hline
$E_3$ & $\rho = 0$   & $\rho = 0.15$                            & $\rho = [-1, 1]$ \\ \hline
$E_4$ & $\rho = 0$   & $\rho = [-0.2, 0.2]$                               & $\rho = [-1, 1]$ \\ \hline
$E_5$ & $\rho = 0$   & $\rho = 1$                               & $\rho = [-1, 1]$ \\ \hline
\end{tabular}
\end{table}

\subsection{Interval scenario}

On the left of Figure 5 it is shown the probability of failure $E_1$ for the interval scenario, where the red area belongs to the independence case, blue to the mixed dependence case, and green to the unknown dependence.
It is possible to see how the uncertainty increases as the dependence assumptions are removed from the Boolean operations.
Intervals for the probability of $E_1$ are indicated in Table 5.

\begin{table}[!th]
\centering
\caption{Probability of $E_1$ in the interval scenario.}\label{tab1}
\begin{tabular}{|c|c|}
\hline
Case               & $\mathbb{P}(E_1)$                   \\ \hline
Independence       & ${[}2.950\times10^{-5}, 4.053\times10^{-5}{]}$ \\ \hline
Mixed Dependence   & ${[}2.949\times10^{-5}, 6.551\times10^{-5}{]}$ \\ \hline
Unknown Dependence & ${[}2.499\times10^{-5}, 1.905\times10^{-4}{]}$  \\ \hline
\end{tabular}
\end{table}

\subsection{P-box scenario}

The results in the P-box scenario are more complicated to interpret since the calculations are probabilities of the probability of the event $E_1$.
However, meaningful analyses can be drawn from them.
For example, assuming the required failure probability of the tank system is no more than $10^{-4}$, a possible question could be what is the probability of it being lower or equal than $10^{-4}$, or how likely is the system to comply with that requirement.
Table 6 includes the results of such inquiry.
In the case of all the events being independent, the probability of fulfilling the requirements goes from 0.969 to 1, so one can infer the tank pressure system is likely to fulfil the requirements.
When the independence assumption is removed, and the mixed dependence scenario is adopted, it is possible to see how the uncertainty on the assessment increases.
Finally, if all the dependence assumptions are relaxed, the probability goes from $0$ to $1$, meaning that no guarantees can be given on the reliability of the tank pressure system.
These results are illustrated in Figure 5 (right), which depicts how dramatic the consequences can be when assuming certain dependencies for the events.

\begin{table}[!th]
\centering
\caption{Probability of $\mathbb{P}(E_1) <= 10^{-4}$ in the P-box scenario.}\label{tab1}
\begin{tabular}{|c|c|}
\hline
Case               & $\mathbb{P}(E_1) <= 10^{-4}$ \\ \hline
Independence       & {[}0.969, 1{]}    \\ \hline
Mixed Dependence   & {[}0.88, 1{]}     \\ \hline
Unknown Dependence & {[}0, 1{]}        \\ \hline
\end{tabular}
\end{table}

%

\begin{figure}[!ht]
    \centering
        \centering
        \includegraphics[width=0.5\linewidth]{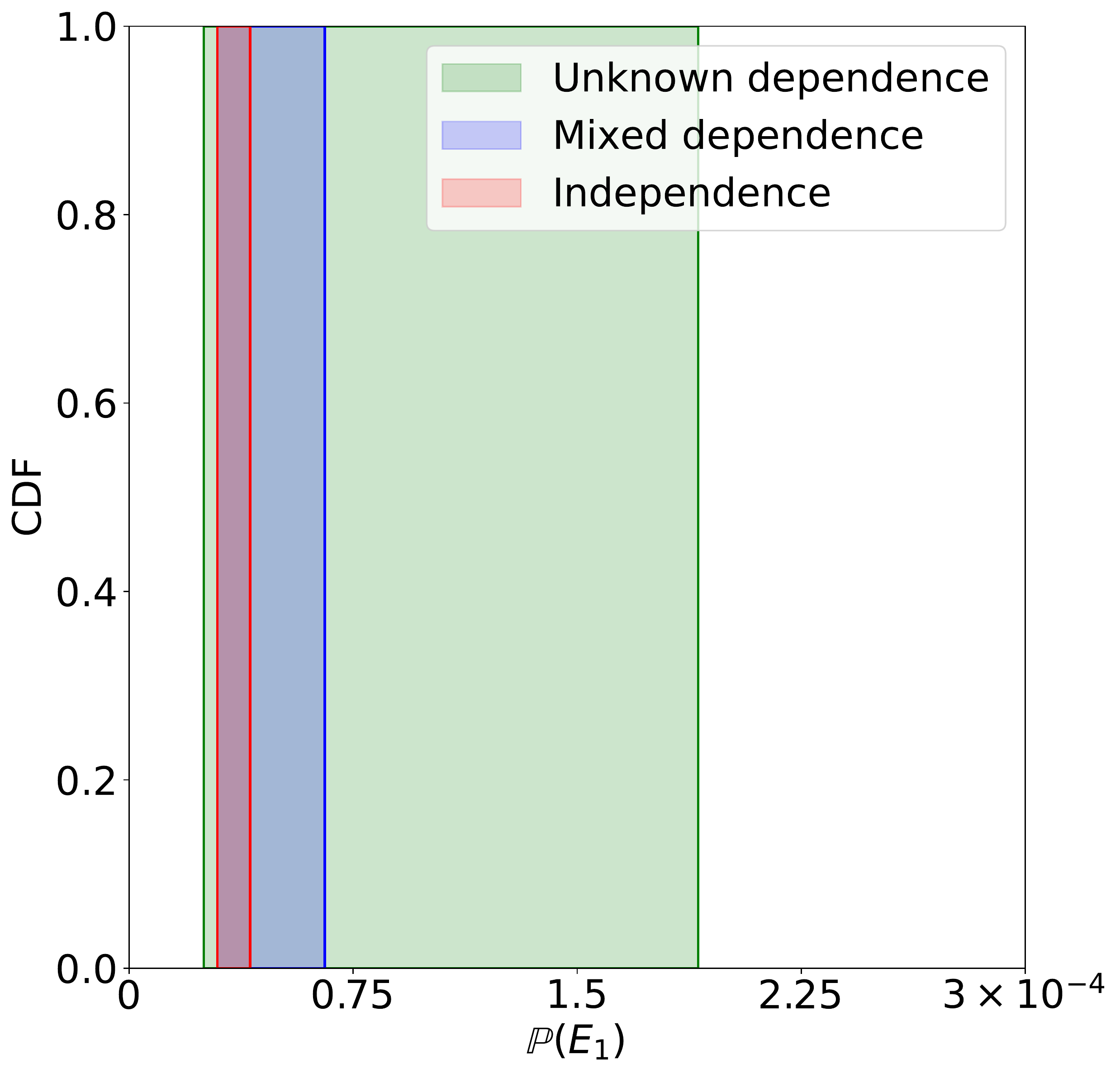}%
        \hfill
        \includegraphics[width=0.5\linewidth]{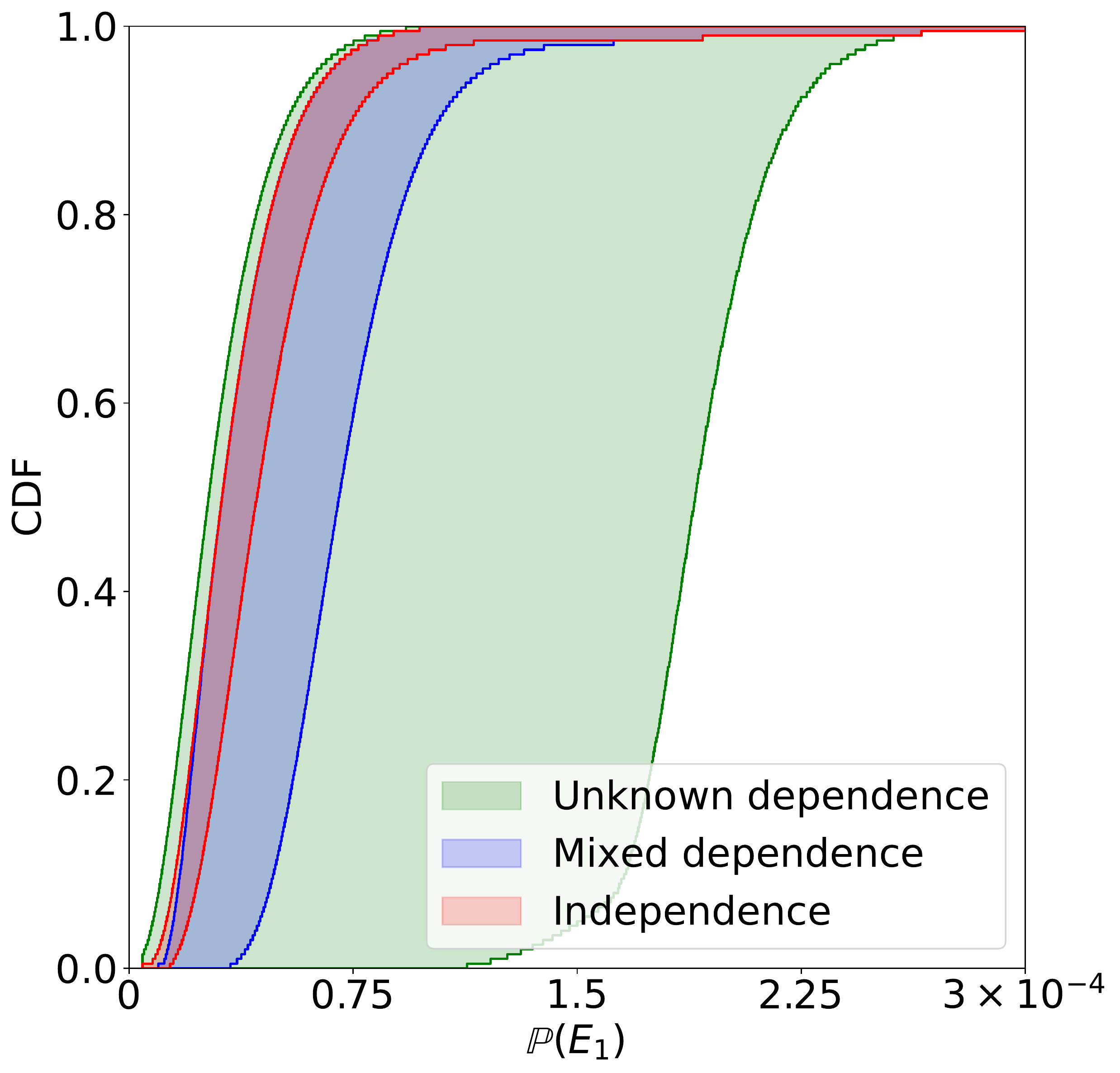}
    \caption{\textbf{(Left)} Interval of $\mathbb{P}(E_1)$ for different dependence assumptions. \textbf{(Right)} Probability box of $\mathbb{P}(E_1)$ for different dependence assumptions.}
\end{figure}

\section{Conclusion}

In this paper, classical Boolean functions have been generalised to be able to operate with precise probabilities, intervals, p-boxes, and with any input correlation.
We show an application of the generalisation of Boolean operations calculating the probability of failure of a pressure tank system.
The uncertainty is propagated through the fault tree with different combinations of probabilities of the events, intervals, and p-boxes, under different dependence assumptions.
The results suggest that assumptions on the probability or dependence of the events have a strong impact on the outcome of the analysis, and these should
be carefully addressed in any serious assessment.

Also, the issue of the two levels of dependence has been introduced, and shown that it can have a dramatic effect on the results of a probabilistic risk assessment.

The computational resources to perform correlated Boolean operations are available in the following open-source Julia package: \url{https://github.com/Institute-for-Risk-and-Uncertainty/UncLogic.jl}

\section*{Acknowledgements}

We thank William (Bill) Huber (from Analysis \& Inference) for his advice in the early stages of this research. This research was funded by the EPSRC and ESRC CDT in Risk and Uncertainty (EP/L015927/1), established within the Institute for Risk and Uncertainty at the University of Liverpool.
This work has been carried out within the framework of the EUROfusion Consortium, funded by the European Union via the Euratom Research and Training Programme (Grant Agreement No 101052200 — EUROfusion). Views and opinions expressed are however those of the author(s) only and do not necessarily reflect those of the European Union or the European Commission. Neither the European Union nor the European Commission can be held responsible for them.

%
%
%

\bibliographystyle{splncs04}
\bibliography{references.bib}

\end{document}